\documentclass[oneside,11pt]{amsart}
\usepackage[english]{babel}
\usepackage{graphicx}

\usepackage{indentfirst} 
\newtheorem{thm}{Theorem}[section]
\newtheorem{cor}[thm]{Corollary}
\newtheorem{lem}[thm]{Lemma}
\newtheorem{prop}[thm]{Proposition}
\theoremstyle{definition}
\newtheorem{defn}[thm]{Definition}
\theoremstyle{remark}
\newtheorem{rem}[thm]{Remark}
\numberwithin{equation}{section}
\begin{document}

\title[Q-P-F property of a class of saddle point matrices]{Quasi-Perron-Frobenius property of a class of saddle point matrices}%
\author{Zheng Li$^{\ast\dag}$, Tie Zhang$^{\dag}$ and Chang-Jun Li$^{\ddag}$}
\thanks{
\indent
$\ast$ Corresponding author, email: neu\_lizheng@hotmail.com,\quad lizheng\_mail@sina.com\\\indent
\dag Department of Mathematics, Northeastern University, Shenyang, 110004, P.R.China.\\\indent
\ddag School of Computer and Software Engineering, University of Science and Technology Liaoning, Anshan, 114051, P.R.China.}


\begin{abstract}
The saddle point matrices arising from many scientific computing fields have block structure
$
W= \left(\begin{array}{cc}
A & B\\
B^T & C
\end{array}
\right)
$,
where the sub-block $A$ is symmetric and positive definite, and $C$ is symmetric and semi-nonnegative definite. In this article we report a unobtrusive but potentially theoretically valuable conclusion that under some conditions, especially when $C$ is a zero matrix, the spectral radius of $W$ must be the maximum eigenvalue of $W$. This characterization approximates to the famous Perron-Frobenius property, and is called quasi-Perron-Frobenius property in this paper. In numerical tests we observe the saddle point matrices derived from some mixed finite element methods for computing the stationary Stokes equation. The numerical results confirm the theoretical analysis, and also indicate that the assumed condition to make the saddle point matrices possess quasi-Perron-Frobenius property is only sufficient rather than necessary.
\end{abstract}

\maketitle
\section{Introduction}\label{sec1}
\subsection{A brief introduction to Perron-Frobenius theory}
Let $\mathbb{R}^n$ and $\mathbb{R}^{m, n}$ denote the spaces of real (column) $n$-vectors and $m\times n$ real matrices respectively, and $\rho(M)$ be the spectral radius of the square matrix $M$. In 1907, Oskar Perron \cite{Perron(1907)} published a fundamental discovery on the positive matrices (whose entries are all real positive numbers). Some results are as follows.

\begin{thm}\label{th0}\cite[part of Th.8.2.8]{Horn&Johnson(2015)}
Let $M\in\mathbb{R}^{n, n}$ be positive. Then\\
(a) $\rho(M)>0$.\\
(b) $\rho(M)$ is an algebraically simple eigenvalue of $M$.\\
(c) There is a unique $x=(x_1,..., x_n)^T\in \mathbb{R}^n$ such that $Mx= \rho(M)x$ and $x_1+...+x_n=1$; this vector is positive.\\
(d) $\lvert \lambda\rvert < \rho(M)$ for every eigenvalue $\lambda$ of $M$ such that $\lambda \not= \rho(M)$.
\end{thm}

Soon after, Georg Frobenius \cite{Frobenius(1912)} generalized Perron' theorem to the case of irreducible and nonnegative matrices (whose entries are all real nonnegative numbers). Subsequently, the famous Perron-Frobenius theory, as an important part of nonnegative matrix theory, was gradually established. Among the results of this theory \cite{Varga(2006), Horn&Johnson(2015), Zhan(2013)}, there is a concise conclusion as follows.

\begin{thm}\label{th0+1}\cite[Th.6.10]{Zhan(2013)}
If $M$ is a nonnegative square matrix, then $\rho(M)$ is an eigenvalue of $M$, and $M$ has a nonnegative eigenvector corresponding to $\rho(M)$.
\end{thm}

Correspondingly, the following classical definition (can be found in \cite{Noutsos(2006), Elhashash&Szyld(2009)}) was derived.
\begin{defn}\label{def0}
Matrix $M\in\mathbb{R}^{n, n}$ is said to possess {\it Perron-Frobenius property} if $M$ satisfies all of following conditions.\\
(a) $\rho(M)$ is an eigenvalue of $M$;\\
(b) There exists a nonnegative eigenvector $x\in\mathbb{R}^n$ such that $Mx= \rho(M)x$.
\end{defn}

In the past 100 years, Perron-Frobenius theory has been greatly developed (see, e.g., \cite{Noutsos(2006), Zhan_report(2007), Elhashash&Szyld(2009), Gluck(2017), Li&Jia(2021)}), and its scope of application has also been expanded. Dimitrios Noutsos \cite{Noutsos(2006)} has given examples that some matrices with negative entries still possess Perron-Frobenius property. This revealed that Perron-Frobenius property is not exclusive to nonnegative matrices. Further more, sometimes we may only focus on the condition (a) of Definition \ref{def0}, which was customarily regarded as Perron--Frobenius property in some literature (see, e.g., \cite{Zhan_report(2007)}), mainly because the information of eigenvalues is the main object in that application. In this article, in order to distinguish this simplified definition from the classical version, we adopt the following expression.

\begin{defn}\label{def1}
Matrix $M\in\mathbb{R}^{n, n}$ is said to possess {\it quasi-Perron-Frobenius property} if $\rho(M)$ is an eigenvalue of $M$.
\end{defn}

In the next subsection, we will show that a class of saddle point matrices, which are well-known in some applications and generally not nonnegative matrices, also possess quasi-Perron-Frobenius property.

\subsection{Background and some classic results on eigenvalue theory of saddle point matrices}
Let the symbol "$T$" represent the transpose of matrix or vector, and symbols "$O$" and "$I$"  denote the zero matrix and identity matrix, respectively.
A {\it saddle point matrix} is the coefficient matrix of the linear system (usually called {\it augmented system, saddle point system} or {\it KKT system}) derived from the saddle point problems occur in many fields, such as computational fluid dynamics, constrained optimization, digital image processing, economics, etc. (see \cite{Benzi&Golub&Liesen(2005)}), and usually has the form as
\begin{equation}\label{eq1}
W=\left(\begin{array}{cc}
A   & B\\
B^T & C
\end{array}
\right)\in\mathbb{R}^{(m+n), (m+n)},\quad usually\; m\geq n,
\end{equation}
where $A\in \mathbb{R}^{m, m}$ is symmetric and positive definite (SPD) (which means all the eigenvalues of $A$ are real positive numbers), and $C\in \mathbb{R}^{n, n}$ is symmetric and semi-negative definite (which means all the eigenvalues of $C$ are located in the interval $(-\infty,\; 0]$ ).
A special and important version of (\ref{eq1}) is that the sub-block $C$ degenerate to the zero matrix, i.e.,
\begin{equation}\label{eq1+0.1}
W=\left(\begin{array}{cc}
A   & B\\
B^T & O
\end{array}
\right),
\end{equation}
which can be generated from many computational models as well (see, e.g., \cite{Benzi&Golub&Liesen(2005), Li&Li&Evans&Zhang(2003), Li&Li&Shao&Nie&Evans(2004), Li&Li&Evans&Zhang(2005), Li&Zhang&Li(2011), Ruiz&Sartenaer&Tannier(2018)}).

In recent decades, researchers have developed a lot of algorithms for solving the augmented system (see, e.g.,  \cite{Benzi&Golub&Liesen(2005), Saad(2009), Cao(2003), Bai(2009), Li&Li&Evans&Zhang(2003), Li&Li&Shao&Nie&Evans(2004), Li&Li&Evans&Zhang(2005)}). And as far as we know, the current research focuses on the preconditioning methods, which essentially relied on the eigenvalue theory of saddle point matrices (see, e.g., \cite{Rusten&Winther(1992), Silvester&Wathen(1994), Perugia&Simoncini(2000), Simoncini&Benzi(2004), Benzi&Golub&Liesen(2005), Benzi&Simoncini(2006), Axelsson&Neytcheva(2006), Huang&Wu&Li(2009), Bai&Ng&Wang(2009), Li&Zhang&Li(2011), Ruiz&Sartenaer&Tannier(2018), Rozloznik(2018)}).

In many applications, the spectral radius of a large scale matrix is an important information, and can usually be computed efficiently by some "simple" algorithm, e.g., the power method. However, even if the eigenvalues are all real numbers, there is still a trouble that the spectral radius may be the positive largest eigenvalue or the opposite of the negative smallest eigenvalue. To accurately identify them generally requires a certain procedure.

Let $\lambda_{max}(\cdot)$, $\lambda_{min}(\cdot)$ and $\lambda_{j}(\cdot)$ denote the maximum, minimum, and $j$th-largest real eigenvalue of matrix respectively, and $\sigma_{max}(\cdot)$, $\sigma_{min}(\cdot)$ and $\sigma_{j}(\cdot)$ denote the maximum, minimum, and $j$th-largest singular values of matrix respectively. It is well known (see e.g., \cite{Benzi&Golub&Liesen(2005)}) that all the eigenvalues of matrix $W$ defined by (\ref{eq1}) satisfy :
\begin{equation}\label{eq1+0.11}
\lambda_{m+n}(W)\leq ...\lambda_{m+1}(W)\leq 0<\lambda_{m}(W)\leq ...\leq \lambda_{1}(W),
\end{equation}
where $\lambda_{m+n}(W)=\lambda_{min}(W)$ and $\lambda_{1}(W)=\lambda_{max}(W)$.
All eigenvalues of sub-blocks $A$ and $C$ and all singular values of sub-block $B$ are assumed as follows,
\begin{equation}\label{eq1+0.12}
0<\lambda_{min}(A)=\lambda_{m}(A)\leq ...\leq \lambda_{1}(A)=\lambda_{max}(A),
\end{equation}
\begin{equation}\label{eq1+0.13}
\lambda_{min}(C)=\lambda_{n}(C)\leq ...\leq \lambda_{1}(C)=\lambda_{max}(C)\leq 0,
\end{equation}
and
\begin{equation}\label{eq1+0.14}
\sigma_{min}(B)=\sigma_{n}(B)\leq ...\leq \sigma_{1}(B)=\sigma_{max}(B).
\end{equation}

For the case that $C=O,\; A=\eta I$, Bernd Fischer et al. \cite{Fischer&Ramage&Silvester&Wathen(1998)} gave some explicit results (has been included in \cite{Benzi&Golub&Liesen(2005)}) as follows.
\begin{thm}\label{thm0.9}\cite[Part 1 of Th.3.8]{Benzi&Golub&Liesen(2005)}
Let matrix $W$ defined as in (\ref{eq1+0.1}). Suppose $A= \eta I$, and $B$ has rank $n-r$, then the $m+n$ eigenvalues of $W$ are given by
\begin{enumerate}
\item zero with multiplicity $r$,
\item $\eta$ with multiplicity $m-n-r$,
\item $\frac{1}{2}\left(\eta\pm\sqrt{4\sigma_k^2+\eta^2}\right)$ for $k=1, ..., n-r$.
\end{enumerate}
\end{thm}
Under the assumption of Theorem \ref{thm0.9}, it follows that
\begin{equation*}
\lambda_{max}(W)=\frac{1}{2}\left(\eta+\sqrt{4\sigma_1^2+\eta^2}\right),\quad \lambda_{min}(W)=\frac{1}{2}\left(\eta-\sqrt{4\sigma_1^2+\eta^2}\right),
\end{equation*}
which implies that $\rho(W)=\lambda_{max}(W)$ (for $\eta>0$), that is, matrix $W$ possesses quasi-Perron-Frobenius property in this case.

However, the situation seems to be complicated if the sub-block $A$ changes to a more general SPD matrix. Generally, we have to turn to seeking some inexact estimations. To our knowledge, among the existing literature, Torgeir Rusten $\&$ Ragnar Winther \cite{Rusten&Winther(1992)} and David Silvester $\&$ Andrew Wathen \cite{Silvester&Wathen(1994)} respectively gave the following interval estimation of the eigenvalues of matrix $W$, which have been included in \cite{Benzi&Golub&Liesen(2005)} or \cite{Rozloznik(2018)}.
\begin{thm}\label{thm1}\cite{Rusten&Winther(1992), Benzi&Golub&Liesen(2005), Rozloznik(2018)}
Let matrix $W$ and its sub blocks $A$ and $B$ are defined as (\ref{eq1+0.1}), and $B$ be of full-column rank. Suppose all eigenvalues of $A$ and all singular values of $B$ are assigned as (\ref{eq1+0.12}) and (\ref{eq1+0.14}) respectively. Then the spectrum of $W$ satisfies
\begin{align*}
& sp(W)\\
&\subset \left[\frac{1}{2}\left(\lambda_m(A)-\sqrt{\lambda_m^2(A)+4\sigma_1^2(B)}\right),\;  \frac{1}{2}\left(\lambda_1(A)-\sqrt{\lambda_1^2(A)+4\sigma_n^2(B)}\right)\right]\\
&\cup \left[\lambda_m(A),\; \frac{1}{2}\left(\lambda_1(A)+\sqrt{\lambda_1^2(A)+4\sigma_1^2(B)}\right)\right].
\end{align*}
\end{thm}
\begin{thm}\label{thm2}\cite{Silvester&Wathen(1994), Rozloznik(2018)}
Let matrix $W$ and its sub-blocks $A$, $B$ and $C$ are defined as (\ref{eq1}), and $B$ be rank-deficient. Suppose all eigenvalues of $A$ and $C$ and all singular values of $B$ are assigned as (\ref{eq1+0.12}), (\ref{eq1+0.13}) and (\ref{eq1+0.14}) respectively, such that the matrix $B^T B-C$ is SPD with the smallest eigenvalues $\lambda_n(B^T B-C)>0$. Then the spectrum of $W$ satisfies
\begin{align*}
sp(W)\subset &\left[\frac{1}{2}\left(\lambda_1(C)+\lambda_m(A)-\sqrt{(\lambda_m(A)-\lambda_1(C))^2+4\sigma_1^2(B)}\right)\right.,\\ &\left.\frac{1}{2}\left(\lambda_1(A)-\sqrt{\lambda_1^2(A)+4\lambda_n^2(C-B^T B)}\right) \right]\\
&\cup \left[\lambda_m(A),\; \frac{1}{2}\left(\lambda_1(A)+\sqrt{\lambda_1^2(A)+4\sigma_1^2(B)}\right)\right].
\end{align*}
\end{thm}

Theorem \ref{thm1} and Theorem \ref{thm2} undoubtedly have important theoretical value. However, we have to note that they are different from Theorem \ref{thm0.9}, since the bounds of those intervals are not exactly coincide with the specific eigenvalues. For example, the upper bound $\frac{1}{2}\left(\lambda_1(A)+\sqrt{\lambda_1^2(A)+4\sigma_1^2(B)}\right)$ of the intervals described in the above theorems is generally not equal to $\lambda_{max}(W)$, although the former is usually a sharp estimate of the latter. The same is true for the other interval bounds. So, it is difficult for us to observe the quasi-Perron-Frobenius property of matrix $W$ based on the above theorems.

In this article, via a short and elementary analysis, we prove that, for the saddle point matrix $W$ defined as (\ref{eq1}), under some conditions, especially when $W$ is defined as (\ref{eq1+0.1}), it is indeed true that $\rho(W)=\lambda_{max}(W)$, i.e., $W$ must possess quasi-Perron-Frobenius property.

\subsection{Organization of the rest of this article}

The remaining sections of this article are organized as follows. In Proposition \ref{prop1} of Section \ref{sec2}, based on Corollary \ref{cor1-1}, we theoretically prove that under some conditions, the saddle point matrices must possess quasi-Perron-Frobenius property. Consequently, in Corollary \ref{cor1} we naturally come to the conclusion that the saddle point matrix defined as (\ref{eq1+0.1}) (i.e. the case $C=O$), must possess quasi-Perron-Frobenius property. In Section \ref{sec3}, via numerical tests we observe the eigenvalue of saddle point matrices generated by some mixed finite element methods for the stationary Stokes equation. The numerical results not only confirm the theoretical conclusion, but also demonstrate that the assumption given in Proposition \ref{prop1} is only sufficient rather than necessary. Finally in Section \ref{sec4}, we briefly summarize the main work of this paper.

\section{Main result}\label{sec2}
\begin{lem}\label{lem1}
Let matrix $W$ defined as (\ref{eq1}). Then\\
(i) $\rho(W)=\lambda_{max}(W)$ if and only if $\lambda_{max}(W)+ \lambda_{min}(W)\geq 0$;\\
(ii) $\rho(W)=\lambda_{max}(W)\not= -\lambda_{min}(W)$ if and only if $\lambda_{max}(W)+ \lambda_{min}(W)> 0$.
\end{lem}
\begin{proof} Since
\begin{align*}
&\rho(W)= \lambda_{max}(W)\Longleftrightarrow \max\left\{\lambda_{max}(W),\; -\lambda_{min}(W) \right\}= \lambda_{max}(W)\\
&\Longleftrightarrow \lambda_{max}(W)\geq -\lambda_{min}(W)\Longleftrightarrow  \lambda_{max}(W)+ \lambda_{min}(W)\geq 0,
\end{align*}
then conclusion (i) is proved. The conclusion (ii) is clearly true.
\end{proof}

Let the inner product between vectors $x$ and $y$ in $\mathbb{R}^n$ be expressed by
$\left<x, y\right>$,
and let
$\|x\|=\sqrt{\left<x, x\right>}$
represent the Euclidean norm of vector. The symbol $\lvert \cdot \rvert$ denotes the absolute value of real number. A useful lemma is as following
\begin{lem}\label{lem2}
Let $B\in \mathbb{R}^{m, n}$. Then
$$\max\limits_{\|x\|^2+\|y\|^2=1}\lvert\left<x, By \right>\rvert= \sigma_{max}(B)/2,$$
and the maximum value is reached when
\begin{equation}\label{eq1+0.15}
\|x\|=\|y\|=\sqrt{2}/2.
\end{equation}
\end{lem}
\begin{proof}
For any real numbers $\alpha>0$, $\beta>0$ we have
\begin{align}
& \max\limits_{\|x\|=\alpha,\;\|y\|=\beta}\lvert \left<x, By \right> \rvert= \max\limits_{\|x/\alpha\|=1,\;\|y/\beta\|=1}\lvert \alpha\beta\left<x/\alpha, By/\beta\right>\rvert \notag\\
& =\alpha\beta\max\limits_{\|x\|=1,\;\|y\|=1}\lvert \left<x, By \right>\rvert =\alpha\beta\sigma_{max}(B).\label{eq1+0.2}
\end{align}
The last equation of (\ref{eq1+0.2}) is derived from \cite[Th.5.6.2 and Page 346]{Horn&Johnson(2015)}.
Therefore, it follows that
\begin{align*}
\max\limits_{\|x\|^2+\|y\|^2=1}\lvert\left<x, By \right> \rvert= \max\limits_{\stackrel{\|x\|^2=\alpha^2,\; \|y\|^2= 1-\alpha^2}{\alpha\in(0, 1)}} \lvert \left<x, By \right> \rvert=\max\limits_{\alpha\in(0, 1)}\alpha \sqrt{1-\alpha^2}\sigma_{max}(B),\notag
\end{align*}
which obviously reaches the maximum when $\alpha=\sqrt{2}/2$, that is,
\begin{equation}\label{eq1+0.3}
\max\limits_{\|x\|^2+\|y\|^2=1}\lvert \left<x, By \right> \rvert= \sigma_{max}(B)/2.
\end{equation}
This complete the proof.
\end{proof}

Consequently, we have the following result.

\begin{cor}\label{cor1-1}
Let $B\in \mathbb{R}^{m, n}$. Then
$$\min\limits_{\|x\|^2+\|y\|^2=1}\left<x, By \right>= -\sigma_{max}(B)/2,$$
and the minimum value is reached when the condition (\ref{eq1+0.15}) holds.
\end{cor}

For the matrix $W$ in (\ref{eq1}), based on above discussion we obtain a characterization of $W$ as follows.
\begin{prop}\label{prop1}
Let matrix $W$ be defined as (\ref{eq1}). If
\begin{equation}\label{eq1+0.5}
\lambda_{min}(A)+\lambda_{min}(C)\geq 0,
\end{equation}
then $W$ possesses quasi-Perron-Frobenius property, i.e.,
\begin{equation*}
\rho(W)=\lambda_{max}(W).
\end{equation*}
In addition, if
\begin{equation*}
\lambda_{min}(A)+\lambda_{min}(C)> 0,
\end{equation*}
then
\begin{equation*}
\rho(W)=\lambda_{max}(W)\not= -\lambda_{min}(W).
 \end{equation*}
\end{prop}
\begin{proof}
For arbitrary vectors $x\in\mathbb{R}^m$, $y\in\mathbb{R}^n$, let $z=(x^T, y^T)^T$, and
\begin{equation*}
F_W(z)= F_W(x, y)=\left<\left(\begin{array}{c} x\\ y \end{array}\right), W \left(\begin{array}{c} x\\ y \end{array}\right)\right>= \left<x, Ax\right>+ 2\left<x, By\right> +\left<y, Cy\right>,
\end{equation*}
then from Courant-Fischer theorem \cite{Horn&Johnson(2015)} it follows that
\begin{equation}\label{eq2}
\lambda_{max}(W)+ \lambda_{min}(W)=\max\limits_{\|x\|^2+\|y\|^2=1}F_W(x, y)+ \min\limits_{\|x\|^2+\|y\|^2=1}F_W(x, y),
\end{equation}
where
\begin{align}\label{eq3}
&\min\limits_{\|x\|^2+\|y\|^2=1}F_W(x, y)\notag \\
&\geq \min\limits_{\|x\|^2+\|y\|^2=1}\left(\left<x, Ax\right> +\left<y, Cy \right> \right)+ \min\limits_{\|x\|^2+\|y\|^2=1}\left<2x, By \right>.
\end{align}
According to Corlloary \ref{cor1-1}, suppose vectors $x_1\in\mathbb{R}^m$ and $y_1\in\mathbb{R}^n$ satisfying
\begin{equation*}
 \|x_1\|^2=1/2,\quad \|y_1\|^2=1/2,
\end{equation*}
are the solutions of the problem
\begin{equation*}
\min\limits_{\|x\|^2+\|y\|^2=1}\left<2x, By \right>.
\end{equation*}
At the same time, there exist vectors $x_2\in\mathbb{R}^m$ and $y_2\in\mathbb{R}^n$, which are the eigenvectors of $\lambda_{min}(A)$ and $\lambda_{min}(C)$ respectively, and satisfy
\begin{equation}\label{eq4}
\|x_2\|= \|x_1\|=\sqrt{2}/2,\quad \|y_2\|= \|y_1\|=\sqrt{2}/2.
\end{equation}
Therefore,
\begin{align}
&\min\limits_{\|x\|^2+\|y\|^2=1}\left(\left<x, Ax\right> +\left<y, Cy \right> \right)+ \min\limits_{\|x\|^2+\|y\|^2=1}\left<2x, By \right>\notag\\
&= \lambda_{min}(A)\|x_2\|^2+ \lambda_{min}(C)\|y_2\|^2+ 2\left<x_1, By_1\right>\notag\\
&= \lambda_{min}(A)\|x_1\|^2+ \lambda_{min}(C)\|y_1\|^2+ 2\left<x_1, By_1\right>.\label{eq5}
\end{align}
Meanwhile, we have
\begin{align}
&\max\limits_{\|x\|^2+\|y\|^2=1}F_W(x, y)\geq F_W(-x_1, y_1) \notag \\
&= \left<x_1, Ax_1\right>- 2\left<x_1, By_1\right> + \left<y_1, Cy_1\right> \notag\\
&\geq \lambda_{min}(A)\|x_1\|^2- 2\left<x_1, By_1\right>+ \lambda_{min}(C)\|y_1\|^2.\label{eq6}
\end{align}
Substituting results of (\ref{eq5}) and (\ref{eq6}) into (\ref{eq2}), and from (\ref{eq4}) we obtain
\begin{align}
\lambda_{max}(W)+ \lambda_{min}(W)&\geq 2\lambda_{min}(A)\|x_1\|^2+ 2\lambda_{min}(C)\|y_1\|^2\notag \\
&= \lambda_{min}(A)+\lambda_{min}(C).\label{eq7}
\end{align}
Combining (\ref{eq7}) with Lemma \ref{lem1}, we complete this proof.
\end{proof}
\begin{rem}\label{remark1}
A further question is: is the sufficient condition (\ref{eq1+0.5}) also necessary for Proposition \ref{prop1}?
In the next section, the numerical results give a negative answer.
\end{rem}

\begin{cor}\label{cor1}
Let matrix $W$ be defined as (\ref{eq1+0.1}). Then $W$ possesses quasi--Perron--Frobenius property, that is,
\begin{equation*}
\rho(W)=\lambda_{max}(W)\not= -\lambda_{min}(W).
\end{equation*}
\end{cor}
\begin{proof}
For the matrix $W$ defined by (\ref{eq1+0.1}), since $C=O$, then the inequality
\begin{equation*}
\lambda_{min}(A)+\lambda_{min}(C)= \lambda_{min}(A)>0
\end{equation*}
always holds. Thus according to Proposition \ref{prop1}, the proof is completed.
\end{proof}
\begin{rem}\label{remark2}
Corollary \ref{cor1} tells us that the spectral radius of matrix $W$ in (\ref{eq1+0.1}) is always the maximum eigenvalue of $W$, other than the absolute value of the minimum eigenvalue of $W$. This may bring convenience in some applications.
\end{rem}
\begin{rem}\label{remark3}
Proposition \ref{prop1} and Corollary \ref{cor1} characterize an exact relation between the spectral radius and the (maximum) eigenvalue of matrix $W$, which is different from Theorem \ref{thm1} and Theorem \ref{thm2}. Moreover, since the sub-block $A$ in Proposition \ref{prop1} or Corollary \ref{cor1} is assumed to be a more general SPD matrix, the corresponding conclusions do not seem to be covered by Theorem \ref{thm0.9}.
\end{rem}

\section{Numerical tests}\label{sec3}

\subsection{Two numerical examples for the stationary Stokes equation}

Consider the stationary Stokes equation \cite{Benzi&Golub&Liesen(2005), Zhang(1992)}
\begin{equation}\label{eq8}
\left\{\begin{array}{l}
-\tau\Delta{\bf u}+\bigtriangledown p = {\bf f},\quad in\; \Omega,\\
div\; {\bf u}=0,\quad in\; \Omega,\\
\int_{\Omega} pd\Omega=0,\\
{\bf u}= {\bf g},\quad on\; \partial \Omega.
\end{array}
\right.,
\end{equation}
where $\Omega=(0, 1)\times(0, 1)$ is a unit square domain, and $\partial \Omega$ is the boundary of $\Omega$. Vector ${\bf u}$ represents the velocity, and $p$ stands for the pressure. The constant $\tau$ is the viscosity coefficient. We apply two mixed finite element methods to discretize the equation (\ref{eq8}). These methods are briefly described as below.

{\bf Example 1: P1-P0 mixed finite element method.}
For the P1-P0 method, we divide $\Omega$ into uniform grids of triangular elements and join the midpoints of the edges on each triangle. We partition each coarse triangle into four refined triangles, and then use the piecewise linear elements for ${\bf u}$ on the fine grid and the piecewise constant elements for $p$ on the coarse grid to discretize the equation (\ref{eq8}) (see \cite{Zhang(1992)} for details). The resulting saddle point matrices have essentially the same form described in (\ref{eq1+0.1}), that is,
\begin{equation}\label{eq9}
W=\left(%
\begin{array}{cc}
  \tau A & B \\
  B^T & O \\
\end{array}%
\right)=\left(%
\begin{array}{ccc}
  \tau \widetilde{A} &  & B_1 \\
   & \tau \widetilde{A} & B_2 \\
  B_1^T & B_2^T & O \\
\end{array}%
\right),
\end{equation}
where $A\in \mathbb{R}^{m\times m}$ is SPD and $B\in \mathbb{R}^{m\times n}$ has full column rank. We use this numerical example to verify Corollary \ref{cor1}. The values of the sum
$S(W)=\lambda_{max}(W)+\lambda_{min}(W)$
under various matrix sizes and viscosity coefficients are reported in Table \ref{tab1}.

{\bf Example 2: Stabilized Q1-P0 mixed finite element method.}
For the stabilized Q1-P0 method, we partition the domain $\Omega$ into $ne\times ne$ uniform square elements, and use piecewise bi-linear elements for ${\bf u}$ and piecewise constant elements for $p$ to discretize the equation (\ref{eq8}). Let $h=1/ne$, the matrix $C$ is resulted from the global stabilization (some details can be found in \cite{Cao(2003)}). Then the resulting coefficient matrices have essentially the same form described in (\ref{eq1}), that is,
\begin{equation}\label{eq10}
W=\left(%
\begin{array}{cc}
  \tau A & B \\
  B^T & C \\
\end{array}%
\right)=\left(%
\begin{array}{ccc}
  \tau\widehat{A} &  & B_1 \\
   & \tau\widehat{A} & B_2 \\
  B_1^T & B_2^T & C \\
\end{array}%
\right).
\end{equation}
where $A\in \mathbb{R}^{m\times m}$ is SPD, $B\in \mathbb{R}^{m\times n}$ has full column rank, and $C\in \mathbb{R}^{n\times n}$ is symmetric and semi-nonnegative definite. This numerical example is used to verify Proposition \ref{prop1}.
In the numerical tests we observe the eigenvalues of the corresponding matrices with different values of viscosity coefficient $\tau$. All the eigenvalue information of these examples are observed on MATLAB. The signs of the sum $S(W)=\lambda_{max}(W)+\lambda_{min}(W)$ and values of $\lambda_{min}(A)$ and $\lambda_{min}(C)$ under various matrix sizes and viscosity coefficients are reported in Table \ref{tab2}.

\subsection{Numerical results}

In Table \ref{tab1} and Table \ref{tab2}, the symbol $m$ represents the order of sub-block $A$ derived in different grids (after deleting some rows and columns from the original generated matrices according to the boundary condition of (\ref{eq8})), while the symbol $n$ represents the number of columns of the sub-block $B$ derived in different subdivision scales (after deleting a column from the original generated matrices according to the condition $\int_{\Omega} pd\Omega=0$). The meaning of other notations and symbols are the same as that in equation (\ref{eq8}) and Proposition \ref{prop1}.

\begin{table}[ht]
\begin{center}
\caption{Numerical results of P1-P0 method}\label{tab1}
\begin{tabular}{c c c c c c}
\hline
{$m$} & $n$ & $\tau$ & $\lambda_{max}(W)$ & $\lambda_{min}(W)$ & $ S(W)$\\
\hline
98	& 31       & 0.01 &	1.05196296	 & -1.01134827	& 0.04061469\\
	&		   & 0.1  &	1.26859883	 & -0.85462817	& 0.41397066\\
	&		   & 1	  & 7.70326732	 & -0.28324435	& 7.42002296\\
	&		   & 10	  & 76.95595469	 & -0.03142832	& 76.92452637\\
\hline
450	& 127	   & 0.01 &	1.02820054	 & -0.98803749	& 0.04016306\\
	&		   & 0.1  &	1.24330881	 & -0.83351936	& 0.40978945\\
	&		   & 1	  & 7.92507966	 & -0.27501021	& 7.65006945\\
	&		   & 10	  & 79.23160501	 & -0.03050241	& 79.20110260\\
\hline
1922 & 511	   & 0.01 &	1.02227690	 & -0.98222993	& 0.04004697\\
	&		   & 0.1  &	1.23699045	 & -0.82827314	& 0.40871731\\
	&		   & 1	  & 7.98122610	 & -0.27297216	& 7.70825394\\
	&		   & 10	  & 79.80743778	 & -0.03027091	& 79.77716688\\
\hline
7938 & 2047    & 0.01 &	1.02079719	 & -0.98077944	& 0.04001776\\
	&		   & 0.1  &	1.23541129	 & -0.82696361	& 0.40844768\\
	&		   & 1	  & 7.99530382	 & -0.27246396	& 7.72283987\\
	&		   & 10	  & 79.95183045	 & -0.03021301	& 79.92161743\\
\hline
\end{tabular}
\footnotetext{Note: $S(W)=\lambda_{max}(W)+\lambda_{min}(W)$}
\end{center}
\end{table}

\begin{table}[ht]
\centering
\caption{Numerical results of stabilized Q1-P0 method}\label{tab2}
\begin{tabular}{c c c c c c c c}\hline
{$m$} & $n$ & $\tau$ & $\lambda_{min}(A)$ & $\lambda_{min}(C)$ & $\lambda_{max}(W)$ & $\lambda_{min}(W)$ & $S(W) $\\
\hline
98 & 63       & 0.01 & 0.002968	 & -0.120235	& 0.233826	& -0.257780	& $-$\\
   & 		  & 0.1  & 0.029676	 & -0.120235	& 0.486079	& -0.167876	& $+$\\
   &		  & 1	 & 0.296756	 & -0.120235	& 3.819638	& -0.121470	& $+$\\
   &		  & 10   & 2.967561	 & -0.120235	& 38.048896	& -0.120354	& $+$\\
\hline
450 & 255     & 1e-3	& 0.000076	 & -0.030950	& 0.118361	& -0.130019	& $-$\\
	&		  & 0.01	& 0.000764	 & -0.030950	& 0.138775	& -0.114891	& $+$\\
	&		  & 0.1     & 0.007637	 & -0.030950	& 0.429370	& -0.050061	& $+$\\
	&		  & 1	    & 0.076367	 & -0.030950	& 3.953115	& -0.031047	& $+$\\
\hline
1922 & 1023   & 1e-4     & 0.000002	 & -0.007794	& 0.060633	& -0.064140	& $-$\\
	 &		  & 1e-3	 & 0.000019	 & -0.007794	& 0.062515	& -0.062434	& $+$\\
	 &		  & 0.01	 & 0.000192	 & -0.007794	& 0.084063	& -0.048097	& $+$\\
	 &		  & 0.1      & 0.001923	 & -0.007794	& 0.408153	& -0.013340	& $+$\\
	 &		  & 1	     & 0.019230	 & -0.007794	& 3.988164	& -0.007800	& $+$\\
\hline
7938 & 4095   & 1e-4     & 4.816e-07 & -0.001952	& 0.030950	& -0.031527	& $-$\\
	 &		  & 1e-3	 & 0.000005	 & -0.001952	& 0.032840	& -0.029820	& $+$\\
	 &		  & 0.01	 & 0.000048	 & -0.001952	& 0.056839	& -0.017847	& $+$\\
	 &		  & 0.1      & 0.000482	 & -0.001952	& 0.402099	& -0.003396	& $+$\\
	 &		  & 1	     & 0.004816	 & -0.001952	& 3.997034	& -0.001952	& $+$\\
\hline
\end{tabular}
\end{table}

The numerical results reported in Table \ref{tab1} illustrate that the values of $S(W)$ are always positive (may be very close to 0 when $\tau\approx 0$, but still remains positive), which is equivalent to $\rho(W)=\lambda_{max}(W)$, for all the cases. This shows that the saddle point matrices with the form (\ref{eq9}) always possess quasi-Perron-Frobenius property, which is consistent with Corollary \ref{cor1}.

The numerical results reported in Table \ref{tab2} show that the saddle point matrices in (\ref{eq10}) possess quasi-Perron-Frobenius property when $\lambda_{min}(A)+\lambda_{min}(C)>0$. This coincides with the conclusion of Proposition \ref{prop1}.

An "unexpected" phenomenon in Table \ref{tab2} is that the saddle point matrices may still possess quasi-Perron-Frobenius property even though $\lambda_{min}(A)+\lambda_{min}(C)<0$, which occur in the lines (m=98, n=63, $\tau$=0.1), (m=450, n=255, $\tau$=0.1, 0.01), (m=1922, n=1023, $\tau$=0.1, 0.01, 1e-3) and (m=7938, n=4095, $\tau$=0.1, 0.01, 1e-3). This indicates that the condition (\ref{eq1+0.5}) is generally only sufficient rather than necessary for the saddle point matrices to possess quasi-Perron-Frobenius property, and makes us have a clearer understanding of the condition of Proposition \ref{prop1}.

\section{Conclusion}\label{sec4}
The main work of this article is to prove that for the saddle point matrix $W$ defined by (\ref{eq1}), if $\lambda_{min}(A)+\lambda_{min}(C)\geq 0$, especially when $C=O$, then $\rho(W)= \lambda_{max}(W)$. This characterize an exact relation between the spectral radius and the maximum eigenvalue of the saddle point matrices. The numerical results also imply that there may exists more precise criteria, which needs further research in the future.

\section*{Acknowledgements}
This work was supported by the State Key Laboratory of Synthetical Automation for Process Industries Fundamental Research Funds, China, no. 2013ZCX02, and National Natural Science Foundation of China, nos: 61575090, 61775169.




\end{document}